\documentclass{amsart}
\usepackage{amsfonts,amssymb,amscd,amsmath,enumerate,verbatim,calc}

\newcommand{\wrt}{with respect to}

\newcommand{\m}{\mathfrak{m} }

 \newcommand{\rt}{\rightarrow}
\newcommand{\xar}{\longrightarrow}

\newcommand{\wh}{\widehat }

\newcommand{\coker}{\operatorname{coker}}

\newcommand{\depth}{\operatorname{depth}}

\newcommand{\Hom}{\operatorname{Hom}}

\theoremstyle{plain}
\newtheorem{thm}{Theorem}
\newtheorem{lem}[thm]{Lemma}

\theoremstyle{definition}

\newtheorem{rem}[thm]{Remark}

\theoremstyle{remark}

\begin{document}

\title [Non-vanishing Theorem  ]{An elementary proof  of  Grothendieck's Non-vanishing Theorem}
\author{Tony ~J.~Puthenpurakal}
\date{\today}
\address{Department of Mathematics, Indian Institute of Technology Bombay, Powai, Mumbai 400 076}

\email{tputhen@math.iitb.ac.in}
 \begin{abstract}
 We give an elementary proof of Grothendieck's non-vanishing Theorem: For a finitely generated non-zero module $M$ over a Noetherian local ring $A$ with maximal ideal
 $\m$, the local cohomology module $H^{\dim M}_{\m}(M)$ is non-zero.
\end{abstract}
\dedicatory{Dedicated to Prof. K. D. Joshi}
 \maketitle
 Let
$(A,\mathfrak{m})$
be a Noetherian local ring, $M$
a finitely generated \textit{non-zero} $A$-module of dimension $r$. For $i \geq 0$ let $H_{\m}^{i}(M)$ be the $i^{th}$ local cohomology module
of $M$ \wrt \ $\m$. Let $\ell(N)$ denote length of an $A$-module $N$.

It is easy to see that $H_{\m}^{i}(M) = 0$ for $i > r$.  Grothendieck's non-vanishing Theorem states that $H^{r}_{\m}(M) \neq 0$.
 However the two well-known proofs of this theorem are quite involved, for instance see \cite[3.5.7(b)]{BH} and \cite[7.3.2]{BSh}.
 The essential point in our proof of the non-vanishing Theorem is that it is easier  to prove by induction the following  stronger result:
\begin{thm}\label{mine}
Let $(A,\mathfrak{m})$
be a Noetherian local ring, $M$
a finitely generated \textit{non-zero} $A$-module of dimension $r$.
  If $r = \dim M \geq 1$ then $ \ell \left(H^{r}_{\m}(M) \right) = \infty $.
\end{thm}
 Theorem \ref{mine} implies Grothendieck's non-vanishing Theorem
since when $\dim M = 0 $ then $H^{0}_{\m}(M) = M$.  Theorem \ref{mine}  is usually deduced as a consequence of the
 non-vanishing Theorem; see \cite[6.6.5]{BSh}.

\begin{rem}\label{pos-depth}
We will need the following well-known basic facts regarding local cohomology:
\begin{enumerate}[\rm (1)]
  \item  $H^{i}_{\m}(M)$ are Artinian for all $i \geq 0$; see \cite[3.5.4(a)]{BH}.

  \item  Let $\wh{M}$ denote the $\m$-adic completion of a finitely generated $A$-module $M$. Then
  \[
  H^{i}_{\m}(M) \cong H^{i}_{\m}(M)\otimes_A \wh{A} \cong H^{i}_{\wh{\m}}(\wh{M}) \quad \text{for all} \ i \geq 0; \ \text{see \cite[3.5.4(d)]{BH}.}
  \]
  \item  When $\dim M > 0$,  set $N = M/H^{0}_{\m}(M)$. Then $\dim N = \dim M$ and $\depth N >0$.
Furthermore $H^i_\m(N) = H^{i}_\m(M)$ for all $i \geq 1$; see \cite[2.1.7]{BSh}.
\item  $H^{i}_{\m}(M) = 0$ for $i< \depth M$ and $i > \dim M$; see \cite[3.5.7(a)]{BH}.
\end{enumerate}
\end{rem}

For $x \in A$ let $\mu^{x}_{M} \colon M \rt M$ be multiplication by $x$. The following lemma is well-known
\begin{lem}\label{mu-f} Let $(A,\mathfrak{m})$
be a Noetherian local ring and let $M$ be
a finitely generated $A$-module. Then
there exists $x \in \m$ such that $\ker \mu^{x}_{M}$ has finite length.
\end{lem}
\begin{proof}\textit{(sketch)}
 If $\dim M = 0$ then any $x \in \m$ will do the job. When $\dim M > 0$ then set $N = M/H^{0}_{\m}(M)$; see  \ref{pos-depth}.3. If $x$ is $N$-regular then it is easy to see
(for instance by the snake lemma) that $\ker \mu^{x}_{M}$ has finite length.
\end{proof}

\begin{proof}[Proof of Theorem \ref{mine}]
By \ref{pos-depth}.2  we may assume $A$ is complete.
 We prove the assertion by induction on $r = \dim M \geq 1$.

\textit{For $r = 1$:} By   \ref{pos-depth}.3 we may assume $\depth M > 0$.
 Let $x \in \m$ be $M$-regular. Set $L = M/xM$. Since $\depth M > 0$ we have $H^{0}_{\m}(M) = 0$. Also as $\dim L = 0$ we have $H^{0}_{\m}(L) =
L$ and $H^{1}_{\m}(L) = 0$.

The short exact sequence
$0 \rt M \xrightarrow{\mu^{x}_{M}} M \rt L \rt 0 $
yields
\[
 0 \xar L \xar  H^{1}_{\m}(M) \rt  H^{1}_{\m}(M) \xar 0.
\]
If $\ell (H^{1}_{\m}(M))$ is finite then  $L =0$. So $M =xM$ and therefore by
Nakayama Lemma we get $M = 0$; a contradiction. Thus $\ell (H^{1}_{\m}(M)) = \infty$.

\textit{We assume the result for modules of dimension $s$ and prove for modules having dimension $s+1$.}
Let $M$ be a module of dimension $s+1$. By  \ref{pos-depth}.3 we may assume $\depth M > 0$.

Let $E$ be the injective hull of the field $A/\m$. By \ref{pos-depth}.1 and  Matlis duality,  the $A$-modules
$H^{i}_{\m}(M)^{\vee} = \Hom_A( H^{i}_{\m}(M), E)$ are finitely generated for all $i \geq 0$.

Set $H = H^{s}_{\m}(M)$ and  $D = M \oplus H^{\vee}$.
By Lemma \ref{mu-f}  there exists
$x \in \m$ such that $\ker \mu^{x}_{D}$ has finite length. It follows that
$\ker \mu^{x}_{H^\vee}$ has finite length. By Matlis duality we get that $C = \coker \mu^{x}_{H}$ has
finite length.

Since $\depth M > 0$ we also get that $x$ is necessarily $M$-regular. Set $L = M/xM$. Notice $L$ is a
module of dimension $s$. Since $s \geq 1$,  by induction hypothesis we get $\ell(H^{s}_{\m}(L)) = \infty$. Also note that
$H^{s+1}_{\m}(L) = 0$.

The short exact sequence
$0 \rt M \xrightarrow{\mu^{x}_{M}} M \rt L \rt 0 $
yields
$$ 0 \rt C \rt H^{s}_{\m}(L) \rt H^{s+1}_{\m}(M)  \rt H^{s+1}_{\m}(M) \rt 0. $$
Since $\ell (C)$ is finite and $\ell (H^{s}_{\m}(L))$ is infinite it follows that
 $\ell(H^{s+1}_{\m}(M))$ is infinite.

Thus by induction we get that if $r = \dim M > 0$ then $\ell(H^{r}_{\m}(M)) = \infty$.
\end{proof}

\begin{rem}
As pointed by the referee, an argument similar to one given above also  proves  Grothendieck's vanishing theorem for finitely generated modules of dimension $r \geq 1$.
\end{rem}
\bigskip\noindent
{\bf Acknowledgment :} The author owes an intellectual debt to Prof K. D. Joshi who taught him that it is sometimes easier to prove stronger results by induction. He also thanks the referee for many pertinent comments.
\providecommand{\bysame}{\leavevmode\hbox to3em{\hrulefill}\thinspace}
\providecommand{\MR}{\relax\ifhmode\unskip\space\fi MR }
\providecommand{\MRhref}[2]{%
  \href{http://www.ams.org/mathscinet-getitem?mr=#1}{#2}
}
\providecommand{\href}[2]{#2}

\end{document}